\documentclass[12pt]{article}

\usepackage[T1]{fontenc}
\usepackage[utf8]{inputenc}
\usepackage{lmodern}
\usepackage{amsmath,amsthm,amsfonts,amssymb,mathtools}
\usepackage{booktabs,array}
\usepackage[a4paper,margin=27mm]{geometry}
\usepackage{microtype}
\usepackage[numbers,sort&compress]{natbib}
\usepackage{needspace,etoolbox}
\usepackage[hidelinks]{hyperref}

\numberwithin{equation}{section}

\newtheorem{teorema}{Theorem}[section]
\newtheorem{proposicion}[teorema]{Proposition}
\newtheorem{lema}[teorema]{Lemma}
\newtheorem{corolario}[teorema]{Corollary}

\theoremstyle{definition}

\newtheorem{ejemplo}[teorema]{Example}
\newtheorem{problema}[teorema]{Problem}
\theoremstyle{remark}
\newtheorem{observacion}[teorema]{Remark}

\AtBeginEnvironment{teorema}{\Needspace{6\baselineskip}}
\AtBeginEnvironment{proposicion}{\Needspace{5\baselineskip}}
\AtBeginEnvironment{lema}{\Needspace{5\baselineskip}}
\AtBeginEnvironment{corolario}{\Needspace{5\baselineskip}}
\AtBeginEnvironment{ejemplo}{\Needspace{5\baselineskip}}
\pretocmd{\subsection}{\Needspace{8\baselineskip}}{}{}
\allowdisplaybreaks[1]

\newcommand{\ind}{\mathbf 1}
\newcommand{\Z}{\mathbb Z}

\newcommand{\N}{\mathbb N}
\newcommand{\U}{\mathsf U}
\newcommand{\D}{\mathcal D}

\newcommand{\ceil}[1]{\left\lceil#1\right\rceil}
\newcommand{\floor}[1]{\left\lfloor#1\right\rfloor}

\hypersetup{
  pdftitle={
    Prime supports and vanishing in
    restricted coprime representations
  },
  pdfauthor={Andres M. Salazar},
  pdfsubject={Paired sieves with fixed prime support},
  pdfkeywords={
    coprime representations,
    prime covers,
    sieve methods,
    Jacobsthal-type bounds,
    residue occupancies
  }
}

\title{
Prime supports and vanishing\\
in restricted coprime representations
}

\author{
Andr\'es M. Salazar\thanks{
\texttt{andresmsalazar@javerianacali.edu.co}
}\\
Department of Natural Sciences and Mathematics\\
Pontificia Universidad Javeriana Cali, Colombia
}

\date{}

\begin{document}
\maketitle
\begin{abstract}
We study the prime supports that eliminate every representation of an even integer $N$ by two positive integers coprime to $6P$. Let $J(k)$ be the largest vanishing instance obtainable with at most $k$ factors of $P$ not dividing $N/2$. We prove
$k(\log k)^2\ll J(k)\ll_\varepsilon k^{4+\varepsilon}$ for every $\varepsilon>0$.
The upper bound follows by comparing a classical family of shifted-prime representations with a two-dimensional upper sieve for the loss at each additional factor. Retaining the same local density factor in both estimates yields a weighted survival criterion and positivity whenever the number of noncolliding factors above $N^\theta$ is sufficiently small compared with $N^\theta$, for fixed $0<\theta<1/4$. An exact reduction separates a minimum prime cover into forced prime-pair contributions and a residual problem, and gives the lower bound for $J$ along powers of two. Residue occupancies characterize the extremal concentration at a first zero. Explicit discrepancy estimates show that supports with at most eleven factors cannot eliminate every representation when $N>2780460$, nor supports with at most seven noncolliding factors when $N\ge2266292$. Positivity for the complete support of primes up to $\sqrt N$ remains open.
\end{abstract}

\section{Introduction}\label{sec:introduccion}
Let $N=2n\ge26$ be even, and let $P$ be a product of distinct primes at least $5$. We consider
\begin{equation}\label{eq:g}
g_P(N)=\#\{1\le h\le n:\gcd(h(N-h),6P)=1\}.
\end{equation}
The endpoint $1+(N-1)$ is included when admissible. The question is how many prime factors, and what arrangement of their excluded classes, are needed to make this count vanish.

The one-prime case was studied in \cite{SalazarCONTACT}. Exact formulas for squarefree supports, together with the parametrization $h=6z+\beta$, discrepancy estimates, and positivity thresholds, were developed in \cite{SalazarCRT}. Here the primes of an instance remain fixed, and the excluded residues are prescribed by that same $N$. Removing a factor reveals the preceding survivor set. Imposing the factors in increasing order identifies a first vanishing step.

Write $\delta_q(x)$ for the least nonnegative remainder of the integer $x$ upon division by $q$. For an odd prime $q\mid P$, the forbidden residues are $0$ and $\delta_q(N)$. They coincide precisely when $q\mid N$, or equivalently $q\mid n$. We call these factors colliding and write
\[
k(N,P)=\omega(P/\gcd(P,n)).
\]
Here $\omega$ counts distinct prime factors. For each integer $k\ge0$, our principal size parameter is
\begin{equation}\label{eq:J}
J(k)=\sup\bigl(\{0\}\cup
\{N\ge26:N\text{ even},\ g_P(N)=0
\text{ for some }P\text{ with }k(N,P)\le k\}\bigr).
\end{equation}
It measures the largest vanishing instance allowed by at most $k$ noncolliding factors, with no bound on their size or on the number of colliding factors. The supremum compares different fixed instances $(N,P)$.

\begin{teorema}[Size of a vanishing instance]\label{thm:main-J}
For every $\varepsilon>0$,
\begin{equation}\label{eq:J-bounds}
k(\log k)^2\ll J(k)\ll_\varepsilon k^{4+\varepsilon}
\qquad(k\longrightarrow\infty).
\end{equation}
In particular, $J(k)$ is finite for every fixed $k$.
\end{teorema}

The upper bound follows from Corollary~\ref{cor:large-factors}. The lower bound is proved in Proposition~\ref{prop:J-lower} by constructing prime covers along the sequence $N=2^j$. This choice avoids a logarithmic loss from the local density factor.

The estimates arise from two complementary descriptions of the same covering problem. First, for any fixed $0<\theta<1/4$, the classical linear sieve supplies
\[
\mathcal T_\theta(N)=\{2\le h<N/2:N-h\text{ prime},\ P^-(h)>N^\theta\},
\]
with cardinality $\gg_\theta\mathfrak S(N)N/\log^2N$, where
$\mathfrak S(N)=\prod_{\ell\mid N,\ \ell>2}(\ell-1)/(\ell-2)$.
An upper sieve applied to $m(N-qm)$ bounds the loss caused by a factor $N^\theta<q\le\sqrt N$ by
$\ll_\theta\mathfrak S(N)N/(q\log^2N)$.
The common local factor cancels. For sufficiently large even $N$, Theorem~\ref{thm:sparse-extensions} gives a weighted positivity criterion and, in particular,
\begin{equation}\label{eq:sparse-intro}
\#\{q\mid P:N^\theta<q<N,\ q\nmid N\}\le\eta_\theta N^\theta
\quad\Longrightarrow\quad
g_P(N)\gg_\theta\mathfrak S(N)\frac{N}{\log^2N}.
\end{equation}
The shifted-prime input is classical \cite{KumchevTolev,BordignonJohnstonStarichkova}. The result here compares its supply with the losses of a fixed support.

Second, let $\tau(N)$ be the minimum total number of prime factors in a support with $g_P(N)=0$. Prime pairs force disjoint covering requirements, each of which can be met by choosing its smaller summand. Theorem~\ref{thm:cover-reduction} separates them exactly:
\begin{equation}\label{eq:cover-intro}
\tau(N)=H(N)+\tau_0(N),
\end{equation}
where $H(N)$ counts the prime pairs with both summands at least $5$, together with the admissible endpoint, and $\tau_0(N)$ is a residual covering number. Fractional covering and disjoint-product packing admit the same reduction, although their optima need not equal the integer covering number. The same prime-pair requirements force $k(N,P)\ge H(N)-1$. Under the conjectural bound $H(N)\gg N/\log^2N$, Remark~\ref{rem:forced-cost} therefore gives $J(k)\asymp k(\log k)^2$.

To describe the first vanishing step, write $P=Qp$ and put
\[
S_Q(N)=\{1\le h<N:\gcd(h(N-h),6Q)=1\}.
\]
Its occupancies satisfy $c_r=c_{\delta_p(N-r)}$. At a first zero, this nonempty set lies entirely in the one or two classes removed by $p$. Theorem~\ref{thm:extremo} characterizes this concentration through a sharp second-moment bound with prescribed central mass. The arithmetic estimates below compare survivor counts with class capacities to exclude such a step.

The distinction between total and noncolliding factors is also useful in explicit ranges. Discrepancy estimates show that supports with at most eleven factors cannot vanish for $N>2780460$, and supports with at most seven noncolliding factors cannot vanish for $N\ge2266292$, with no restriction on the size of the factors. These effective conclusions complement \eqref{eq:J-bounds}, whose constants have not been made explicit.

The parameter $J(k)$ is related in purpose to Jacobsthal-type covering problems. Iwaniec \cite{IwaniecJacobsthal} proved a bound of order $(k\log k)^2$ for the maximal run covered by $k$ primes in the ordinary one-class problem. Ziller and Morack \cite{ZillerMorackTheory} study paired progressions and their phases. Nguyen \cite{NguyenWindows} studies finite-window CRT intersections and shift correlations for symmetric pairs on primorial wheels. Here the classes are tied to the fixed sum $N$, and colliding factors are not charged. The parameter $J(k)$ measures the extremal size under this noncolliding budget.

We call $p_{\max}^2<N$ the natural cutoff. With
$\Pi(N)=\prod_{5\le q\le\sqrt N}q$, factorization gives
\begin{equation}\label{eq:goldbach-endpoint}
g_{\Pi(N)}(N)=\ind_{N-1\text{ prime}}+
\#\{\sqrt N<\lambda\le\varrho:\lambda,\varrho\text{ prime},\ \lambda+\varrho=N\}.
\end{equation}
Every composite summand below $N$ has a prime divisor below $\sqrt N$, so the only admissible summands are $1$ and primes exceeding $\sqrt N$. Since $N$ is even, $\sqrt N$ cannot be an odd prime. Positivity at the complete cutoff therefore implies that every even $N\ge26$ with $N-1$ composite is a sum of two primes exceeding $\sqrt N$. This conclusion remains open here.

Sections~\ref{sec:preliminares}--\ref{sec:momentos} describe the count and the structure of a first zero. Section~\ref{sec:positividad} gives explicit exclusions. Section~\ref{sec:uniformidad} proves the weighted survival bound and the upper bound for $J$. Section~\ref{sec:covers} proves the cover reduction and the lower bound for $J$. Section~\ref{sec:alcance} states the remaining quantitative gap.
\section{Parametrization and discrepancy}\label{sec:preliminares}
Throughout, $N=2n\geq26$. All prime products are squarefree and supported on primes at least $5$. The empty product $Q=1$ is allowed. For $P\ne1$, let $p_{\max}$ be its largest factor. For integers $q\ge1$ and $x\in\Z$, we write the canonical remainder operator as in \cite{SalazarCONTACT,SalazarCRT}:
\begin{equation}\label{eq:delta}
\delta_q(x)=x-q\floor{x/q},\qquad
\mathcal R_q=\{0,1,\ldots,q-1\}.
\end{equation}
Thus $\delta_q:\Z\to\mathcal R_q$, including for negative $x$, and
$\delta_q(x+y)=\delta_q(\delta_q(x)+\delta_q(y))$.
Divisibility is equivalently expressed by $q\mid x$ or $\delta_q(x)=0$. We reserve $\delta$ for this operator and use $\varepsilon>0$ in exponent bounds.

Alongside \eqref{eq:g}, consider the ordered count
\[
s_Q(N)=\#S_Q(N),\qquad D_Q(N)=\ind_{\gcd(n,6Q)=1}.
\]
The reflection $h\mapsto N-h$ has only one possible fixed point, namely $h=n$. Hence
\begin{equation}\label{eq:ordenado}
s_Q(N)=2g_Q(N)-D_Q(N),\qquad
s_Q(N)>0\ \Longleftrightarrow\ g_Q(N)>0.
\end{equation}

Put
\begin{equation}\label{eq:parametros}
\ell=\begin{cases}
n/3,&\delta_3(n)=0,\\
(n+2)/3,&\delta_3(n)=1,\\
(n-2)/3,&\delta_3(n)=2,
\end{cases}
\qquad
\beta=\begin{cases}1,&\delta_3(n)\ne2,\\5,&\delta_3(n)=2,\end{cases}
\qquad b_6=1+\ind_{3\mid N}.
\end{equation}
The parameter interval is $0\leq z<\ell$, with $h=6z+\beta$. If $3\mid n$, this family represents one of the two admissible values of $\delta_6(h)$, and its reflection represents the other. If $3\nmid n$, it represents all permitted summands. Consequently, setting
\begin{equation}\label{eq:CQ}
C_Q(N)=\#\{0\leq z<\ell:\gcd((6z+\beta)(N-6z-\beta),Q)=1\},
\end{equation}
we have
\begin{equation}\label{eq:s-C}
s_Q(N)=b_6C_Q(N),\qquad \ell\geq\frac{N-4}{6}.
\end{equation}
This is the parametrization of \cite{SalazarCRT}, written for the ordered count.

For each $q\mid Q$, let $\iota_q\in\{1,\ldots,q-1\}$ be the unique integer satisfying $\delta_q(6\iota_q)=1$. The excluded parameter residues are
\[
\delta_q(-\iota_q\beta),\qquad \delta_q(\iota_q(N-\beta)).
\]
Their number is $\nu_q=2-\ind_{q\mid N}$. Define
\begin{equation}\label{eq:rhoD}
\rho_Q=\prod_{q\mid Q}\left(1-\frac{\nu_q}{q}\right),\qquad
\D_Q=\prod_{q\mid Q}(1+\nu_q)-\prod_{q\mid Q}\left(1+\frac{\nu_q}{q}\right).
\end{equation}
The dependence on $N$ is suppressed in $\rho_Q$ and $\D_Q$. The discrepancy bound $\D_Q$ is distinct from the diagonal indicator $D_Q(N)$. For $Q=1$, we have $\rho_Q=1$ and $\D_Q=0$.

\begin{lema}[Parameter discrepancy]\label{lem:discrepancia}
We have
\begin{equation}\label{eq:discrepancia}
|C_Q(N)-\ell\rho_Q|\leq\D_Q.
\end{equation}
If $Q\ne1$, then $\D_Q<\prod_{q\mid Q}(1+\nu_q)-1$. In particular,
\begin{equation}\label{eq:lower-s}
s_Q(N)\geq b_6(\ell\rho_Q-\D_Q).
\end{equation}
\end{lema}
\begin{proof}
We recall the proof of \cite[Proposition 5.2]{SalazarCRT}. For each $d\mid Q$, choosing one excluded class at every prime dividing $d$ gives $\prod_{q\mid d}\nu_q$ distinct prescribed values of $\delta_d(z)$ by the Chinese remainder theorem. The count of each prescribed value in an interval of length $\ell$ differs from $\ell/d$ by at most $1-1/d$. Inclusion--exclusion gives the main term
\[
\ell\sum_{d\mid Q}\mu(d)\frac{\prod_{q\mid d}\nu_q}{d}=\ell\rho_Q.
\]
The sum of the remainder bounds is
\[
\sum_{d\mid Q}\left(\prod_{q\mid d}\nu_q\right)\left(1-\frac1d\right)=\D_Q.
\]
The second product in \eqref{eq:rhoD} exceeds $1$ if $Q\ne1$. Finally, apply \eqref{eq:s-C}.
\end{proof}
\section{The structure of a first zero}\label{sec:momentos}
If the factors are imposed in increasing order, the first vanishing step must concentrate the preceding set in the classes removed by the last prime. We first analyze that concentration without an ordering assumption. Let $P=Qp$, where $p\nmid6Q$. We do not yet require $p$ to be the largest factor of $P$. Abbreviate
\begin{equation}\label{eq:ocupaciones}
S=S_Q(N),\qquad s=|S|,\qquad
c_r=\#\{h\in S:\delta_p(h)=r\},\quad r\in\mathcal R_p.
\end{equation}
Reflection of the summands gives
\begin{equation}\label{eq:simetria}
c_r=c_{\delta_p(N-r)}.
\end{equation}

All occupancy indices belong to $\mathcal R_p$. For each phase $a\in\mathcal R_p$, define
\begin{equation}\label{eq:fases}
Y(a)=\begin{cases}
s-c_{\delta_p(-a)}-c_{\delta_p(N-a)},&p\nmid N,\\
s-c_{\delta_p(-a)},&p\mid N.
\end{cases}
\end{equation}
These phases arise by translating $h$ by multiples of $6Q$. Every condition at $Q$ remains unchanged, while $\delta_p(6Qj)$ runs through $\mathcal R_p$ as $j$ does. Translated values need not lie in the original interval. The phase counting the representations of $N$ is $a=0$.

Put
\begin{equation}\label{eq:momentos}
\nu=2-\ind_{p\mid N},\quad
\mu_p=s\left(1-\frac\nu p\right),\quad
V_p=\frac1p\sum_{a\in\mathcal R_p}(Y(a)-\mu_p)^2,
\end{equation}
and
\begin{equation}\label{eq:AB}
A=\sum_{r\in\mathcal R_p}c_r^2,\qquad B=\sum_{r\in\mathcal R_p}c_rc_{\delta_p(r+N)}.
\end{equation}
The mean of $Y$ is $\mu_p$, since each occupancy occurs $\nu$ times in the sum of the losses.

\begin{proposicion}[Variance identity and pointwise bound]\label{prop:varianza}
Let $t=Y(0)=s_P(N)=2g_P(N)-D_P(N)$. Then
\begin{equation}\label{eq:V}
V_p=\begin{cases}
\displaystyle\frac2p(A+B)-\frac{4s^2}{p^2},&p\nmid N,\\[4pt]
\displaystyle\frac Ap-\frac{s^2}{p^2},&p\mid N,
\end{cases}
\end{equation}
and
\begin{equation}\label{eq:cota-t}
t\geq\begin{cases}
\displaystyle s-\sqrt{\frac{2p}{3}V_p+\frac{8s^2}{3p}},&p\nmid N,\\[5pt]
\displaystyle s-\sqrt{pV_p+\frac{s^2}{p}},&p\mid N.
\end{cases}
\end{equation}
\end{proposicion}
\begin{proof}
Expanding the square of the losses in \eqref{eq:fases} gives \eqref{eq:V}. If $p\nmid N$, symmetry yields $c_0=c_{\delta_p(N)}=(s-t)/2$. Nonnegativity gives
\[
A\geq\frac{(s-t)^2}{2},\qquad B\geq\frac{(s-t)^2}{4}.
\]
Substitution into \eqref{eq:V} yields
$V_p+4s^2/p^2\geq3(s-t)^2/(2p)$, proving the first bound. If $p\mid N$, then $c_0=s-t$ and $A\geq(s-t)^2$, which proves the second.
\end{proof}

For $s>0$, write $\U(N,Q,p)$ for the inequality
\begin{equation}\label{eq:U}
V_p<\begin{cases}
\displaystyle\frac{3p-8}{2p^2}s^2,&p\nmid N,\\[4pt]
\displaystyle\frac{p-1}{p^2}s^2,&p\mid N.
\end{cases}
\end{equation}
By the preceding proposition, $\U(N,Q,p)$ implies $g_P(N)>0$. Using \eqref{eq:V}, its two branches become
\begin{equation}\label{eq:energia}
A+B<\frac34s^2\quad(p\nmid N),\qquad A<s^2\quad(p\mid N).
\end{equation}

\subsection{Equality structure}
When the two excluded classes coincide, the criterion has an immediate interpretation.

\begin{proposicion}[Collision]\label{prop:colision}
If $p\mid N$ and $s>0$, then $A\leq s^2$ and
\[
\U(N,Q,p)\quad\Longleftrightarrow\quad g_P(N)>0.
\]
\end{proposicion}
\begin{proof}
The identity $s^2-A=2\sum_{r<r'}c_rc_{r'}$ shows that equality holds precisely when a single class carries mass. By \eqref{eq:simetria}, this class must satisfy $\delta_p(2r)=0$. Since $p$ is odd and $r\in\mathcal R_p$, it is $0$. This is equivalent to $t=0$.
\end{proof}

For the other branch, write $b=\delta_p(N)\ne0$ and
\begin{equation}\label{eq:central}
u=c_{\delta_p(n)},\qquad v=s-u.
\end{equation}
The residue $\delta_p(n)$ is the unique fixed point of $r\mapsto\delta_p(N-r)$. The quantity $u$ counts an entire class, possibly containing several integers $n+jp$. It need not equal the diagonal indicator $D_Q(N)$.

\begin{teorema}[Second moment with prescribed central mass]\label{thm:extremo}
Let $p\geq5$ be prime, let $b\in\mathcal R_p\setminus\{0\}$, and let the nonnegative real numbers $c_r$ satisfy $c_r=c_{\delta_p(b-r)}$. Put
\[
r_*=\delta_p\bigl((p+1)b/2\bigr),\qquad
s=\sum_{r\in\mathcal R_p}c_r,\qquad u=c_{r_*},\qquad v=s-u.
\]
Then
\begin{equation}\label{eq:extremo}
\mathcal E(c):=\sum_r c_r^2+\sum_r c_rc_{\delta_p(r+b)}
\leq\max\left\{u^2+\frac34v^2,\ u^2+uv+\frac12v^2\right\}.
\end{equation}
The bound is attained for every pair of real numbers $u,v\geq0$. If $u=0$ and $s>0$, then
\begin{equation}\label{eq:igualdad}
\mathcal E(c)\leq\frac34s^2,
\end{equation}
and equality holds precisely when $c_0=c_b=s/2$ and all other occupancies vanish.
\end{teorema}
\begin{proof}
The quadratic form can be written as
\begin{equation}\label{eq:convexidad}
\mathcal E(c)=\frac12\sum_{r\in\mathcal R_p}(c_r+c_{\delta_p(r+b)})^2.
\end{equation}
It is positive definite. Indeed, a vector in its kernel would satisfy $c_{\delta_p(r+b)}=-c_r$ for every $r$. Addition by $b$ traverses a cycle of odd length $p$, so $c_r=-c_r$ and the vector vanishes. Thus $\mathcal E$ is strictly convex.

The reflection $r\mapsto\delta_p(b-r)$ has one fixed point $r_*$ and $(p-1)/2$ two-point orbits. Once $u$ is fixed, the occupancies are described by distributing the mass $v$ among these orbits, with equal mass at the two points of each orbit. The set of distributions is a simplex. Convexity allows $\mathcal E$ to be maximized at a vertex, with all the mass $v$ on a single orbit.

For the orbit $\{0,b\}$, the two points are adjacent under the step $b$, and neither is adjacent to the fixed point. The value is $u^2+3v^2/4$. For the orbit $\{\delta_p(-r_*),\delta_p(3r_*)\}$, both points are adjacent to the fixed point but not to each other, and the value is $u^2+uv+v^2/2$. Every other orbit gives $u^2+v^2/2$. These classifications use $p\geq5$. The two special orbits exist and are distinct even for $p=5$. This proves the bound and its sharpness.

If $u=0$, the only vertex attaining $3s^2/4$ is the orbit $\{0,b\}$. Strict convexity rules out equality at a nontrivial combination of vertices, proving the equality characterization.
\end{proof}

\begin{corolario}[Central survival or strict variance]\label{cor:central-or}
Suppose $p\nmid N$ and $s>0$. Then
\[
g_{Qp}(N)>0\quad\Longleftrightarrow\quad
u>0\ \text{or}\ \U(N,Q,p),
\qquad u=c_{\delta_p(n)}.
\]
\end{corolario}
\begin{proof}
The central residue is distinct from $0$ and $\delta_p(N)$, so $u>0$ supplies a survivor. If $u=0$, Theorem~\ref{thm:extremo} shows that failure of the strict criterion is equivalent to concentration in the excluded classes. The converse implications follow from these same observations and Proposition~\ref{prop:varianza}.
\end{proof}

For example, $N=26$, $Q=35$, and $p=11$ give $S_Q(N)=\{13\}$. Here $s=u=1$, so $A+B=1>3/4$, although $g_{385}(26)=1$. The variance criterion is therefore sufficient in general and exact in the empty-central-class branch. The positivity estimates below are obtained independently from the complete survivor set.
\section{Class capacity and explicit positivity}\label{sec:positividad}
Each value of $\delta_p(h)$ is taken by a limited number of integers $h$ in the interval with $\gcd(h(N-h),6)=1$. This allows us to compare a lower bound for the preceding set with an upper bound for the loss.

\begin{lema}[Residue-class capacity]\label{lem:capacidad}
With the notation of \eqref{eq:ocupaciones}, let
\begin{equation}\label{eq:M}
M_{N,p}=b_6\ceil{\frac{N-1}{6p}}.
\end{equation}
Then
\begin{equation}\label{eq:capacidad}
\max_r c_r\leq M_{N,p},\qquad A\leq M_{N,p}s,\qquad B\leq M_{N,p}s.
\end{equation}
\end{lema}
\begin{proof}
The condition $\gcd(h(N-h),6)=1$ permits $b_6$ values of $\delta_6(h)$. Together with a fixed value of $\delta_p(h)$, each determines one value of $\delta_{6p}(h)$. An interval of length $N-1$ contains at most $\lceil(N-1)/(6p)\rceil$ integers with a prescribed value of $\delta_{6p}(h)$. The restrictions at $Q$ can only decrease the count. The bounds for $A$ and $B$ follow by bounding one factor of each summand by $M_{N,p}$.
\end{proof}

\begin{proposicion}[Two capacity criteria]\label{prop:capacidad}
If $s>0$, the conditions
\[
s>\frac83M_{N,p}\quad(p\nmid N),\qquad
s>M_{N,p}\quad(p\mid N)
\]
are sufficient for $\U(N,Q,p)$. To obtain $g_{Qp}(N)>0$, it is enough that
\begin{equation}\label{eq:cap-positividad}
s>\nu M_{N,p},\qquad \nu=2-\ind_{p\mid N}.
\end{equation}
\end{proposicion}
\begin{proof}
The first assertions follow from \eqref{eq:capacidad} and \eqref{eq:energia}. For the last one, symmetry gives
\begin{equation}\label{eq:perdida}
t=s-\nu c_0\geq s-\nu M_{N,p}.
\end{equation}
\end{proof}

Thus an estimate based on the maximum occupancy can prove positivity under a weaker hypothesis than that needed for the variance criterion. The independent bound \eqref{eq:lower-s} allows both comparisons.

\subsection{Complete blocks of the preceding modulus}
\begin{proposicion}\label{prop:bloques}
Suppose $Q>1$ and $p$ exceeds every prime factor of $Q$. Let
$A_Q=Q\rho_Q=\prod_{q\mid Q}(q-\nu_q)$. Then
\begin{equation}\label{eq:bloques}
C_{Qp}(N)\geq\floor{\frac\ell Q}
\left[A_Q-\nu\ceil{\frac Qp}\right],
\qquad \nu=2-\ind_{p\mid N}.
\end{equation}
The bracket is at least $1$. In particular, $\ell\geq Q$ implies
$g_{Qp}(N)>0$, and the same holds if $N\geq6Q$.
\end{proposicion}
\begin{proof}
Each block of $Q$ consecutive parameters contains exactly $A_Q$ survivors of the restrictions at $Q$. Each value of $\delta_p(z)$ occurs at most $\lceil Q/p\rceil$ times in the block. Summing the complete blocks proves \eqref{eq:bloques}.

Let $A_0=\prod_{q\mid Q}(q-2)$. This is an odd integer at least $3$, and $A_Q\geq A_0$. The telescoping product over odd integers from $5$ to $p-2$ gives
\[
\frac{A_0}{Q}\geq\frac3{p-2},\qquad
\frac Qp<\frac{A_0}{3}\leq\frac{A_0-1}{2}.
\]
Since the last quantity is an integer, $2\lceil Q/p\rceil\leq A_0-1$. The bracket in \eqref{eq:bloques} is therefore at least $1$. Finally, the three branches of \eqref{eq:parametros} show that $N\geq6Q$ implies $\ell\geq Q$.
\end{proof}

\subsection{Thresholds for a given support}
Assume now that every prime factor of $Q$ is less than $p$. Put
\begin{equation}\label{eq:rho0}
\rho_0=\prod_{q\mid Q}\left(1-\frac2q\right),\qquad
\D_0=3^{\omega(Q)}-\prod_{q\mid Q}\left(1+\frac2q\right).
\end{equation}
We have $\rho_Q\geq\rho_0$ and $\D_Q\leq\D_0$. For the second comparison, the expression for $\D_Q$ in \eqref{eq:rhoD} increases with each $\nu_q$. Also,
\begin{equation}\label{eq:telescopica}
\rho_0\geq\prod_{\substack{5\leq j\leq p-2\\j\text{ odd}}}\frac{j-2}{j}
=\frac3{p-2}.
\end{equation}
The product on the right includes additional factors smaller than one.

\begin{teorema}[Explicit thresholds]\label{thm:umbrales}
Let $Q$ be a product of primes in $[5,p)$, with $p\ge5$ prime, and let $\rho_0,\mathcal D_0$ be as in \eqref{eq:rho0}. Then
\begin{equation}\label{eq:Tpos}
N>\frac{4\rho_0+6\D_0+12-2/p}{\rho_0-2/p}
\quad\Longrightarrow\quad g_{Qp}(N)>0.
\end{equation}
Furthermore,
\begin{equation}\label{eq:Tvar}
N>\frac{4\rho_0+6\D_0+16-8/(3p)}{\rho_0-8/(3p)}
\quad\Longrightarrow\quad s>0\text{ and }\U(N,Q,p).
\end{equation}
Both denominators are positive.
\end{teorema}
\begin{proof}
By Lemma~\ref{lem:discrepancia},
$s\geq b_6((N-4)\rho_0/6-\D_0)$. Also,
$M_{N,p}\leq b_6((N-1)/(6p)+1)$. Comparing the first expression with $2M_{N,p}$ gives \eqref{eq:Tpos}, and comparing it with $8M_{N,p}/3$ gives \eqref{eq:Tvar}. The latter comparison also covers the branch $p\mid N$. Finally, \eqref{eq:telescopica} implies $\rho_0>8/(3p)>2/p$.
\end{proof}

The term $\D_0$ retains exponential growth in the number of prime factors of $Q$. The theorem uses the fact that the last prime need not be included in the inclusion--exclusion remainder. This improvement of the finite bounds does not by itself cover all supports at the natural cutoff.

\subsection{Supports with at most eleven prime factors}
Let $\ell_1=5,\ell_2=7,\ldots$ be the primes starting at $5$. For $r\geq1$, define
\begin{equation}\label{eq:Tr}
\rho_r=\prod_{i<r}\left(1-\frac2{\ell_i}\right),\quad
E_r=3^{r-1}-1,\quad
T_r=\frac{4\rho_r+6E_r+12-2/\ell_r}{\rho_r-2/\ell_r}.
\end{equation}

\begin{lema}\label{lem:cola-r}
If $\omega(P)=r$ and $N>T_r$, then $g_P(N)>0$.
\end{lema}
\begin{proof}
Remove the largest prime $p$. The prime factors of $Q=P/p$ are at least $\ell_1,\ldots,\ell_{r-1}$, and $p\geq\ell_r$. Hence $\rho_Q\geq\rho_r$ and $\D_Q\leq E_r$. Using \eqref{eq:lower-s} and \eqref{eq:perdida},
\[
t\geq\frac{b_6}{6}\left[
N\left(\rho_r-\frac2{\ell_r}\right)-4\rho_r-6E_r-12+\frac2{\ell_r}
\right].
\]
The substitution $p\geq\ell_r$ is valid in the term $-2(N-1)/p$. The right-hand side is positive for $N>T_r$.
\end{proof}

Table~\ref{tab:colas} lists the ceilings of the first eleven bounds. All comparisons use exact fractions. In particular,
\begin{equation}\label{eq:Tonce}
\max_{1\leq r\leq11}T_r=T_{11}
=\frac{94655100382586}{34042961}
=2780460.265562\ldots.
\end{equation}

\begin{table}[htbp]
\centering\small
\begin{tabular}{rrrrrrrrrrrr}
\toprule
$r$&1&2&3&4&5&6&7&8&9&10&11\\
\midrule
$\ell_r$&5&7&11&13&17&19&23&29&31&37&41\\
$\lceil T_r\rceil$&26&84&250&860&2754&9359&29745&90603&292506&893241&2780461\\
\bottomrule
\end{tabular}
\caption{Thresholds by cardinality. The sufficient condition is $N>T_r$.}
\label{tab:colas}
\end{table}

\begin{teorema}[Eleven prime factors]\label{thm:once}
If $1\leq\omega(P)\leq11$ and $N>2780460$, then $g_P(N)>0$. No restriction on the size of the factors is needed.
\end{teorema}
\begin{proof}
By \eqref{eq:Tonce}, every even $N>2780460$ satisfies $N>T_r$ for $1\le r\le11$. Apply Lemma~\ref{lem:cola-r}.
\end{proof}

\subsection{Prime factors dividing half the integer}
Discrepancy also allows us to retain information about collisions. Write
\begin{equation}\label{eq:perfil}
C=\gcd(P,n),\quad R=P/C,\quad c=\omega(C),\quad k=\omega(R).
\end{equation}
The prime factors of $C$ and $R$ form disjoint sets, $C\mid n$, and $\gcd(R,n)=1$. For this support, Lemma~\ref{lem:discrepancia} implies
\begin{equation}\label{eq:perfil-criterio}
\ell\prod_{q\mid C}\left(1-\frac1q\right)
\prod_{q\mid R}\left(1-\frac2q\right)\geq2^c3^k-1
\quad\Longrightarrow\quad g_P(N)>0.
\end{equation}
The non-strict hypothesis is valid because the discrepancy bound by $2^c3^k-1$ is strict.

For $c\geq1$, let $C_c=\prod_{i=1}^c\ell_i$, and define
\begin{equation}\label{eq:Bck}
B_{c,k}=\frac{\varphi(C_c)}3\left(1-\frac2{C_c}\right)
\prod_{i=c+1}^{c+k}\left(1-\frac2{\ell_i}\right).
\end{equation}

\begin{lema}[Extrema with disjoint groups]\label{lem:grupos}
If $c\geq1$, the left-hand side of \eqref{eq:perfil-criterio} is at least $B_{c,k}$. Moreover, for all $c,k\geq0$ with $c+k\geq1$, a sufficient condition for positivity is
\begin{equation}\label{eq:Theta}
n\geq\Theta_{c,k}:=
3(2^c3^k-1)
\prod_{i=1}^k\frac{\ell_i}{\ell_i-2}
\prod_{i=k+1}^{k+c}\frac{\ell_i}{\ell_i-1}+2.
\end{equation}
\end{lema}
\begin{proof}
Since $n\geq C$ and $\ell\geq(n-2)/3$, the left-hand side of \eqref{eq:perfil-criterio} is at least
\begin{equation}\label{eq:Bactual}
\frac{\varphi(C)}3\left(1-\frac2C\right)
\prod_{q\mid R}\left(1-\frac2q\right).
\end{equation}
Replacing a prime by a smaller unused prime decreases this expression, regardless of its group. Its minimum therefore uses the first $c+k$ primes. If $a<b$, assigning $a$ to $C$ and $b$ to $R$ decreases the product of factors, since
\[
(a-1)\frac{b-2}{b}<(b-1)\frac{a-2}{a}.
\]
It also decreases $C$, and hence the factor $1-2/C$. The minimizing assignment puts the $c$ smallest primes in $C$, proving the first assertion.

Without using $n\geq C$, criterion \eqref{eq:perfil-criterio} gives the threshold
\[
3(2^c3^k-1)\prod_{q\mid C}\frac q{q-1}
\prod_{q\mid R}\frac q{q-2}+2.
\]
Its factors decrease as the primes increase. If $a<b$, the ratio of the assignment $a\in R,b\in C$ to the opposite assignment is
\[
\frac{(a-1)(b-2)}{(a-2)(b-1)}>1.
\]
The maximum assigns the $k$ smallest primes to $R$ and gives \eqref{eq:Theta}.
\end{proof}

\begin{teorema}[Seven noncolliding factors]\label{thm:siete}
If $\omega(P/\gcd(P,n))\leq7$ and $n\ge1133146$, then $g_P(N)>0$, with no restriction on $\omega(P)$ or on the size of its factors. If at least six factors of $P$ divide $n$, the restriction on $n$ is unnecessary.
\end{teorema}
\begin{proof}
If $P=1$, positivity follows from \eqref{eq:s-C}. Suppose $P>1$. For $c=6,k=7$, we have exactly
\[
B_{6,7}=\frac{34212642163660800}{189431751553}
>139967=2^6 3^7-1.
\]
Decreasing $k$ increases $B_{c,k}$ and decreases $2^c3^k-1$. As $c$ increases, the ratio $B_{c+1,k}/B_{c,k}$ exceeds $\ell_{c+1}-1\geq4$, whereas $2^c3^k-1$ grows by a factor at most $3$ for $c\geq1$. Thus Lemma~\ref{lem:grupos} covers all cases $c\geq6$, $k\leq7$.

If $c\leq5$ and $k\leq7$, the maximum $\Theta_{c,k}$ does not exceed $\Theta_{5,7}$. Monotonicity follows by adding factors to the maximization problem in the lemma. The exact value is
\[
\Theta_{5,7}=\frac{23314807798777987}{20575296000}<1133146.
\]
Consequently, all $n\geq1133146$ are covered.
\end{proof}

\begin{teorema}[Restrictions on any zero]\label{thm:restricciones}
If $n\geq13$ and $g_P(2n)=0$, then $k\geq1$ and
\begin{align}
(c+1)!&<\frac{5(2k+3)}3\,3^k,\label{eq:factorial}\\
n&<(2(c+k)+3)(2^c3^k-1)+2.\label{eq:cardinal-n}
\end{align}
These conditions do not require the natural cutoff.
\end{teorema}
\begin{proof}
If $k=0$, every prime factor of $P$ divides $n$. The displacement
$d=2^{\ind_{2\nmid n}}3^{\ind_{3\nmid n}}\leq6<n$
produces the pair $n-d,n+d$, coprime to $6P$. This contradicts vanishing.

For $c\geq1$, the $i$th prime factor of $C$ is at least $2i+3$, so
\[
\varphi(C)\geq2^c(c+1)!,\qquad 1-2/C\geq3/5.
\]
Furthermore,
\[
\prod_{q\mid R}\left(1-\frac2q\right)
\geq\prod_{i=1}^k\frac{2i+1}{2i+3}=\frac3{2k+3}.
\]
Expression \eqref{eq:Bactual} is therefore at least
$3\,2^c(c+1)!/(5(2k+3))$. If \eqref{eq:factorial} failed, the threshold in \eqref{eq:perfil-criterio} would be reached. For $c=0$, \eqref{eq:factorial} is immediate because $k\geq1$.

Put $r=c+k$ and $\rho=\prod_{q\mid C}(1-1/q)\prod_{q\mid R}(1-2/q)$. Vanishing requires $(n-2)\rho/3<2^c3^k-1$. Moreover,
\[
\rho^{-1}\leq\prod_{q\mid P}\frac q{q-2}
\leq\prod_{i=1}^r\frac{2i+3}{2i+1}=\frac{2r+3}{3}.
\]
Combining these bounds proves \eqref{eq:cardinal-n}.
\end{proof}

The elementary restrictions above are effective and retain information about both groups of primes. Equation~\eqref{eq:profile-asymptotic} below gives a stronger asymptotic restriction by using the distribution of primes.
\section{Sparse supports and the size of a vanishing instance}\label{sec:uniformidad}
We now estimate the number of prime factors needed to eliminate every representation. The argument constructs a family with a prime summand and compares its size with the loss at each factor. Keeping the same local density factor in both estimates is essential for uniformity.

\subsection{A family supplied by the linear sieve}
For $h>1$, let $P^-(h)$ denote its least prime factor. Fix $0<\theta<1/4$ and put
\begin{equation}\label{eq:prime-family}
\mathcal T_\theta(N)=\{2\le h<N/2:N-h\text{ prime},\ P^-(h)>N^\theta\}.
\end{equation}
Write
\begin{equation}\label{eq:singular}
\mathfrak S(N)=\prod_{\substack{\ell\mid N\\\ell>2}}\frac{\ell-1}{\ell-2},
\qquad
\mathfrak S_w(N)=\prod_{\substack{\ell\mid N\\2<\ell<w}}\frac{\ell-1}{\ell-2},
\end{equation}
where the products run over primes. We have $\mathfrak S(N)\ge1$.

\begin{lema}[Shifted-prime sieve]\label{lem:shifted-primes}
For every $0<\theta<1/4$, there are $c_\theta>0$ and $N_0(\theta)$ such that
\begin{equation}\label{eq:prime-mass}
|\mathcal T_\theta(N)|\ge c_\theta\mathfrak S(N)\frac{N}{\log^2N}
\end{equation}
for all even $N\ge N_0(\theta)$. Both summands of every representation in this family are coprime to $6N$.
\end{lema}
\begin{proof}
Consider $\mathcal A_N=\{N-r:N/2<r<N,\ r\text{ prime}\}$ and
$X_N=\operatorname{li}(N)-\operatorname{li}(N/2)\asymp N/\log N$.
A prime $r>N/2$ cannot divide the even integer $N$, so every element is coprime to $N$. For squarefree $d$, the divisibility density is
\[
g_N(d)=\begin{cases}1/\varphi(d),&(d,N)=1,\\0,&(d,N)>1.\end{cases}
\]
Thus $|(\mathcal A_N)_d|=X_Ng_N(d)+R_d$. Bombieri--Vinogradov gives, for any $A>0$ and a suitable $B$,
\begin{equation}\label{eq:BV-remainder}
\sum_{d\le D}\mu^2(d)|R_d|\ll_A N\log^{-A}N,
\qquad D=N^{1/2}/\log^B N.
\end{equation}
For $d\le D$ with $(d,N)>1$, the count and its density are both zero. For the other $d$, the estimate is uniform in $N$, since it uses the maximum over reduced classes at both endpoints of the prime interval.

The lower linear sieve applies uniformly after deleting the odd primes dividing $N$. Its density product is
\[
V_N(z)=\prod_{\substack{3\le\ell<z\\\ell\nmid N}}\left(1-\frac1{\ell-1}\right)
\asymp\frac{\mathfrak S_z(N)}{\log z}.
\]
For $z=2N^\theta$,
\[
0\le\log\frac{\mathfrak S(N)}{\mathfrak S_z(N)}
\ll\frac{\log N}{z\log z}=o_\theta(1).
\]
Indeed, at most $\log N/\log z$ prime factors of $N$ are at least $z$.
The sieve parameter satisfies $\log D/\log z\to1/(2\theta)>2$, so its lower function is positive. Consequently,
\[
S(\mathcal A_N,z)\gg_\theta
\mathfrak S(N)\frac{N}{\log^2N}.
\]
Removing the possible element $1$ proves \eqref{eq:prime-mass}. For $N^\theta>3$, the surviving summands are also coprime to $6$.
\end{proof}

This is the classical Goldbach-sequence argument in \cite[Theorem 8 and Sections 4.2--4.3]{KumchevTolev}, with a variable power cutoff. An explicit general form of the linear sieve appears in \cite[Theorem 6]{BordignonJohnstonStarichkova}. The restriction $\theta<1/4$ comes from its lower-sieve parameter and is not a new exponent for that classical construction.

\subsection{The loss at an additional prime}
The roughness of the complementary factor improves the loss estimate by one logarithm. The local factor $\mathfrak S(N)$ must be retained.

\begin{lema}[Two linear forms]\label{lem:upper-forms}
Fix $0<\theta<1/4$ and let $\sigma=\min(\theta/2,1/100)$.
For $q=1$, or for any prime $N^\theta<q\le\sqrt N$ with $q\nmid N$, put $X=N/(2q)$ and $w=N^\sigma$. Then, uniformly in $q$ and for even $N$ sufficiently large,
\begin{equation}\label{eq:upper-forms}
\#\{1\le m<X:\gcd(m(N-qm),\prod_{\ell<w}\ell)=1\}
\ll_\theta\mathfrak S(N)\frac{N}{q\log^2N}.
\end{equation}
\end{lema}
\begin{proof}
No sieving prime divides $q$. The polynomial $m(N-qm)$ has one root at $\ell=2$, one root at each odd $\ell\mid N$, and two roots at all other primes $\ell<w$. Denote this number by $\omega_N(\ell)$ and extend it multiplicatively to squarefree integers. The count for a squarefree $d$ supported below $w$ is
\[
X\frac{\omega_N(d)}d+R_d,\qquad |R_d|\le\omega_N(d)\le2^{\omega(d)}.
\]
The upper-bound sieve of dimension two applies with a bounded upper sieve function, uniformly in these local densities. Indeed, its inverse density products are bounded by those with two roots at every odd prime. Take $D=X/\log^{20}N$. For sufficiently large $N$, one has $40\le\log D/\log w\le1/\sigma$, and
\[
\sum_{d<D}\mu^2(d)4^{\omega(d)}|R_d|
\le\sum_{d<D}8^{\omega(d)}
\ll D(1+\log D)^7\ll X\log^{-13}N.
\]
The standard upper sieve therefore bounds the left side of \eqref{eq:upper-forms} by
\[
\ll_\theta X\prod_{\ell<w}\left(1-\frac{\omega_N(\ell)}\ell\right)
+O(X\log^{-13}N)
\ll_\theta X\frac{\mathfrak S_w(N)}{\log^2w}.
\]
Here the root at $2$ contributes $1/2$, and the correction at every odd divisor of $N$ is $(\ell-1)/(\ell-2)$. Since $\mathfrak S_w(N)\le\mathfrak S(N)$ and $\log w=\sigma\log N$, the claim follows. The upper sieve used here is the standard dimension-two specialization of \cite[Theorem 9.1]{DiamondHalberstam}, with the divisibility remainders verified above.
\end{proof}

Let $L_q(N)=\#\{h\in\mathcal T_\theta(N):q\mid h(N-h)\}$. For a support $P$, define the normalized loss bound
\begin{equation}\label{eq:weighted-cost}
\begin{split}
W_\theta(N,P)={}&
\sum_{\substack{q\mid P\\N^\theta<q\le\sqrt N\\q\nmid N}}\frac1q\\
&+\frac{\log^2N}{\mathfrak S(N)N}
\sum_{\substack{q\mid P\\\sqrt N<q<N\\q\nmid N}}
\left(\frac{N}{2q}+1\right).
\end{split}
\end{equation}
All sums over factors of a support run over primes. Also put
\begin{equation}\label{eq:large-profile}
\kappa_\theta(N,P)=\#\{q\mid P:N^\theta<q<N,\ q\nmid N\}.
\end{equation}

\begin{lema}[Loss at one factor]\label{lem:prime-loss}
For even $N$ sufficiently large,
\[
L_q(N)\le
\begin{cases}
C_\theta\mathfrak S(N)N/(q\log^2N),&N^\theta<q\le\sqrt N,\ q\nmid N,\\
N/(2q)+1,&\sqrt N<q<N,\ q\nmid N,\\
0,&\text{otherwise}.
\end{cases}
\]
In particular, $\max_qL_q(N)\ll_\theta\mathfrak S(N)N^{1-\theta}/\log^2N$.
\end{lema}
\begin{proof}
The zero-loss cases follow from Lemma~\ref{lem:shifted-primes} and the defining roughness. If $N^\theta<q\le\sqrt N$, the prime summand exceeds $N/2$ and cannot equal $q$. Write $h=qm$. Every counted $m$ avoids all prime divisors below $w$ in Lemma~\ref{lem:upper-forms}, so that lemma applies. For $q>\sqrt N$, multiples of $q$ among $h<N/2$ contribute at most $N/(2q)$, and $N-h=q$ contributes at most one. Finally,
$\sqrt N/2+1=o_\theta(N^{1-\theta}/\log^2N)$.
\end{proof}

\begin{teorema}[Weighted survival and sparse supports]\label{thm:sparse-extensions}
Fix $0<\theta<1/4$ and write $x_+=\max\{x,0\}$. There are positive constants $c_\theta,C_\theta,\eta_\theta$ such that, for all sufficiently large even $N$ and every support $P$,
\begin{equation}\label{eq:weighted-survival}
g_P(N)\ge
\mathfrak S(N)\frac{N}{\log^2N}
\bigl(c_\theta-C_\theta W_\theta(N,P)\bigr)_+.
\end{equation}
In particular,
\begin{equation}\label{eq:sparse-hypothesis}
\kappa_\theta(N,P)\le\eta_\theta N^\theta
\quad\Longrightarrow\quad
g_P(N)\gg_\theta\mathfrak S(N)\frac{N}{\log^2N}.
\end{equation}
There is no upper bound on the prime factors of $P$.
\end{teorema}
\begin{proof}
At least $|\mathcal T_\theta(N)|-\sum_{q\mid P}L_q(N)$ members survive. Lemmas~\ref{lem:shifted-primes} and \ref{lem:prime-loss}, enlarging $C_\theta$ to be at least $1$, give \eqref{eq:weighted-survival}. Each nonzero weight in \eqref{eq:weighted-cost} is at most $N^{-\theta}$ once $N$ is sufficiently large. For the second sum this follows from
\[
\frac{\log^2N}{\mathfrak S(N)N}\left(\frac{\sqrt N}{2}+1\right)
=o_\theta(N^{-\theta}).
\]
Hence $W_\theta(N,P)\le\kappa_\theta(N,P)N^{-\theta}$.
Choose $\eta_\theta=c_\theta/(2C_\theta)$.
\end{proof}

The same factor $\mathfrak S(N)$ multiplies the supply and the upper bound for each small-prime loss. Its cancellation explains the absence of a logarithm in \eqref{eq:sparse-hypothesis}. The weighted statement also permits many more factors when they lie close to $\sqrt N$.

\begin{corolario}[Cost and size of a zero]\label{cor:large-factors}
For every $0<\theta<1/4$, a vanishing instance with sufficiently large $N$ satisfies
\begin{equation}\label{eq:zero-large}
\kappa_\theta(N,P)\gg_\theta N^\theta.
\end{equation}
Consequently, for each $\varepsilon>0$ and $k=\omega(P/\gcd(P,N/2))\ge1$,
\begin{equation}\label{eq:zero-polynomial}
g_P(N)=0\quad\Longrightarrow\quad N\ll_\varepsilon k^{4+\varepsilon}.
\end{equation}
Thus $J(k)\ll_\varepsilon k^{4+\varepsilon}$.
\end{corolario}
\begin{proof}
Vanishing requires every member of $\mathcal T_\theta(N)$ to be removed. Compare its lower bound with the maximum loss in Lemma~\ref{lem:prime-loss} multiplied by $\kappa_\theta$. The common factor $\mathfrak S(N)/\log^2N$ cancels. Since $\kappa_\theta\le k$, take $\theta=1/(4+\varepsilon)$. The bounded remaining values of $N$ are absorbed into the implied constant.
\end{proof}

For a sequence of zeros with $N\to\infty$, one also has $k\to\infty$ and
\begin{equation}\label{eq:profile-asymptotic}
c=O\left(\frac{\log k}{\log\log k}\right),
\qquad c=\omega(\gcd(P,N/2)).
\end{equation}
Indeed, $(c+1)!\le\gcd(P,N/2)\le N/2$, whereas \eqref{eq:zero-polynomial} gives $\log N=O(\log k)$.

\subsection{Consequences for the first-zero structure}
The arithmetic estimates also control the localized moments. Their role here is a consequence of the survivor mass, independently of the extremal characterization.

\begin{lema}[Mass for long intervals]\label{lem:long-mass}
Fix $a>4$. There is $p_1(a)$ such that, for every prime $p\ge p_1(a)$, every product $Q$ of primes in $[5,p)$, and every even $N\ge p^a$,
\begin{equation}\label{eq:masa-poly}
s_Q(N)\gg_a\frac{N}{\log^2p}.
\end{equation}
\end{lema}
\begin{proof}
Let $U\ge5$ be an absolute constant, chosen below. If $p^a\le N<p^U$, put $\theta=1/a<1/4$. For $p$ sufficiently large in terms of $a$, Lemma~\ref{lem:shifted-primes} applies. Every $h\in\mathcal T_\theta(N)$ satisfies $P^-(h)>N^\theta\ge p$, and $N-h$ is a prime exceeding $N/2>p$. Thus both orientations survive $Q$, and
\[
s_Q(N)\ge2|\mathcal T_\theta(N)|
\gg_a\frac{N}{\log^2N}
\gg_a\frac{N}{\log^2p},
\]
using $\mathfrak S(N)\ge1$ and $\log N<U\log p$.

Now let $N\ge p^U$. Sieve the parameters $0\le z<\ell$ by the $\nu_q$ excluded values at each $q\mid Q$. Their densities satisfy $\nu_q/q\le2/q$. The inverse-product condition for a sieve of dimension two holds uniformly by comparison with two excluded values at every prime $q\ge5$. For squarefree $d\mid Q$, the divisibility remainder obeys
\[
|R_d|\le\prod_{q\mid d}\nu_q\le2^{\omega(d)}.
\]
Take $D=\ell^{1/2}$ and $s_0=\log D/\log p$. For sufficiently large $p$, one has $s_0\ge U/3$. The lower-bound part of \cite[Theorem~9.1]{DiamondHalberstam}, with lower sieve function $f_2$, gives
\[
C_Q(N)\ge\ell\rho_Q\bigl(f_2(s_0)-o(1)\bigr)
-O\left(\sum_{\substack{d<D\\d\mid Q}}4^{\omega(d)}|R_d|\right).
\]
The sieve error $o(1)$ is $O((\log\log D)^2/(\log D)^{1/6})$, with an absolute implied constant because the inverse-product condition is uniform. Since $f_2(s)$ increases to $1$, choose the absolute $U$ so that $f_2(U/3)>0$. The remainder is bounded by
\[
\sum_{d<D}8^{\omega(d)}\ll D(1+\log D)^7
\ll\ell^{1/2}(1+\log\ell)^7.
\]
Also,
\[
\rho_Q\ge\prod_{5\le q<p}\left(1-\frac2q\right)
\gg\frac1{\log^2p}.
\]
After division by $\ell/\log^2p$, the remainder is
$O(\ell^{-1/2}(1+\log\ell)^9)=o(1)$ uniformly for $N\ge p^U$ as $p\to\infty$.
Hence $C_Q(N)\gg\ell/\log^2p$, and \eqref{eq:s-C} proves the claim.
\end{proof}

\begin{corolario}[Uniformity of increasing prefixes]\label{cor:prefix-uniformity}
Under \eqref{eq:sparse-hypothesis}, $\U(N,Q,p)$ holds at every increasing prefix $Qp$ of $P$ for all sufficiently large $N$, uniformly in $P$. Moreover, for every fixed $a>4$,
\begin{equation}\label{eq:uniforme-poly}
\frac{pV_p}{s_Q(N)^2}\ll_a\frac{\log^2p}{p}
\end{equation}
whenever $N\ge p^a$, $p$ is sufficiently large in terms of $a$, and $Q$ is supported on primes in $[5,p)$.
\end{corolario}
\begin{proof}
Under \eqref{eq:sparse-hypothesis}, the surviving prime family gives
$s_Q(N)\gg_\theta N/\log^2N$ for every $Q\mid P$. Capacity yields
\[
\frac{pV_p}{s_Q(N)^2}\le\frac{4M_{N,p}}{s_Q(N)}
\ll_\theta\log^2N\left(\frac1p+\frac1N\right).
\]
This proves the strict criterion for $p\ge K_\theta\log^2N$, with a sufficiently large fixed $K_\theta$. For the smaller steps, $Q$ is supported below $p$ and $N\ge p^5$ uniformly when $N$ is large. Lemma~\ref{lem:long-mass} with $a=5$ gives $s_Q(N)\gg N/\log^2p$, so the normalized variance is $O(\log^2p/p)$ and satisfies the strict criterion for all sufficiently large $p$. The remaining bounded primes and their finitely many preceding supports are covered by \eqref{eq:Tvar}.

For the second assertion, Lemma~\ref{lem:long-mass} and capacity give
\[
\frac{pV_p}{s_Q(N)^2}\le\frac{4M_{N,p}}{s_Q(N)}
\ll_a\frac{\log^2p}{p}.
\]
Since the right side tends to zero, $\U$ holds for sufficiently large $p$.
\end{proof}

\begin{observacion}[Almost all $N$]\label{rem:almost-all}
For every $A>0$, there is a set of $O_A(X/\log^A X)$ even integers in $[26,X]$ outside which $s_Q(N)>0$ and $\U(N,Q,p)$ hold simultaneously for every prime $5\le p<\sqrt N$ and every $Q$ supported in $[5,p)$. To see this, apply the classical weighted abundance estimate in \cite[proof of Theorem 4.3]{BhowmikGrimmelt} on dyadic intervals. Proper prime powers and primes at most $\sqrt N$ contribute only $O(\sqrt N\log^3N)=o(N)$ to the weighted count. Dividing the remaining weight by $\log^2N$ gives $\gg N/\log^2N$ prime pairs with both summands exceeding $\sqrt N$. These pairs survive every natural-cutoff support. Capacity handles $p>N^{1/5}$, and \eqref{eq:uniforme-poly} handles the smaller unbounded primes. The finitely many remaining primes are covered by \eqref{eq:Tvar}. The exceptional set is therefore common to all supports for the same $N$.
\end{observacion}
\section{Minimum prime covers and a residual problem}\label{sec:covers}
The preceding section bounds from below the noncolliding cost of a zero. We now reduce the minimum total covering problem by selecting one prime from each prime pair. This will give the lower bound for $J(k)$ along powers of two.

Put
\[
E_N=\{1\le h\le N/2:\gcd(h(N-h),6)=1\},\qquad
\mathcal P_h=\{q\ge5:q\text{ prime},\ q\mid h(N-h)\}.
\]
Every $\mathcal P_h$ is nonempty. Define
\begin{equation}\label{eq:cover-number}
\tau(N)=\min\{\omega(P):g_P(N)=0\}.
\end{equation}
This minimum is finite: taking all primes in $[5,N)$ covers every $\mathcal P_h$. Equivalently, $\tau(N)$ is the minimum size of a set of primes meeting every $\mathcal P_h$. Primes at least $N$ can be omitted.

Let $G(N)=g_{\Pi(N)}(N)$ be the count in \eqref{eq:goldbach-endpoint}, and put
\begin{align}
M_0(N)&=\pi(\sqrt N)-2,\label{eq:cover-m}\\
H(N)&=\ind_{N-1\text{ prime}}+
\#\{5\le h\le N/2:h,N-h\text{ prime}\},\label{eq:cover-H}\\
H_{\rm sm}(N)&=\#\{5\le q\le\sqrt N:q,N-q\text{ prime}\}.
\label{eq:cover-a}
\end{align}
Then $H=G+H_{\rm sm}$. Distinct prime pairs use disjoint sets of primes, because each prime determines its complement. The endpoint contributes another disjoint singleton when $N-1$ is prime.

Define the following set of primes:
\[
\mathcal F_N=\{h:5\le h\le N/2,\ h,N-h\text{ prime}\}
\cup\bigl(\{N-1\}\text{ if }N-1\text{ is prime}\bigr).
\]
Thus $|\mathcal F_N|=H(N)$. The residual rows are
\[
E_N^0=\{h\in E_N:\mathcal P_h\cap\mathcal F_N=\varnothing\}.
\]
Let $\tau_0(N)$ be their minimum covering number, with value zero if $E_N^0$ is empty.

\begin{teorema}[Canonical prime-pair reduction]\label{thm:cover-reduction}
A minimum cover can be chosen to contain $\mathcal F_N$, and
\begin{equation}\label{eq:cover-exact}
\tau(N)=H(N)+\tau_0(N),\qquad 0\le\tau_0(N)\le M_0(N)-H_{\rm sm}(N).
\end{equation}
In particular,
\begin{equation}\label{eq:cover-sandwich}
G(N)+H_{\rm sm}(N)\le\tau(N)\le G(N)+M_0(N).
\end{equation}
\end{teorema}
\begin{proof}
Consider a prime pair $5\le\lambda\le\varrho$ with $\lambda+\varrho=N$. If $\lambda<\varrho$, then $\varrho>N/2$. The prime $\varrho$ occurs in only one row of the covering problem, namely the row $h=\lambda$: its only positive multiple below $N$ is $\varrho$. A cover containing $\varrho$ but not $\lambda$ can therefore replace $\varrho$ by $\lambda$ without increasing its size or uncovering any row. If $\lambda=\varrho$, the row forces the prime $\lambda$. The endpoint, when present, forces $N-1$.

The replacements for distinct prime pairs do not interfere. Hence a minimum cover can be chosen to contain $\mathcal F_N$. Its additional primes must cover exactly $E_N^0$, proving the identity. Every residual row has a prime factor at most $\sqrt N$, since otherwise it would be a surviving prime pair or the endpoint already covered by $\mathcal F_N$. None of the $H_{\rm sm}(N)$ small primes in $\mathcal F_N$ occurs in a residual row. The other $M_0(N)-H_{\rm sm}(N)$ small primes cover all residual rows, proving the upper bound.
\end{proof}

The set $\mathcal F_N$ specifies a canonical choice within some minimum cover. The theorem does not require every minimum cover to contain it.

In particular, $\tau(N)>M_0(N)$ implies $G(N)>0$. The converse has not been established here. The unconditional relation $0\le\tau(N)-G(N)\le M_0(N)$ explains why $\tau$ is close to the large-prime representation count whenever that count is large. It does not imply that $G+M_0-\tau$ is bounded, or that all minimum covers have the same shape.

\subsection{Fractional covers and disjoint products}
Let $\nu(N)$ be the maximum cardinality of a family of rows whose products $h(N-h)$ are pairwise coprime. The fractional covering number is
\begin{align}
\tau^*(N)&=\min\left\{\sum_qx_q:x_q\ge0,\ \sum_{q\in \mathcal P_h}x_q\ge1\ \text{for every }h\right\}\label{eq:fractional-cover}\\
&=\max\left\{\sum_hy_h:y_h\ge0,\ \sum_{h:q\in \mathcal P_h}y_h\le1\ \text{for every }q\right\}.
\label{eq:fractional-dual}
\end{align}
The equality is finite-dimensional linear programming duality. Integer covers and disjoint products satisfy
\begin{equation}\label{eq:packing-chain}
H(N)\le\nu(N)\le\tau^*(N)\le\tau(N).
\end{equation}
The duality does not assert equality of the two integer optima.

\begin{proposicion}\label{prop:cover-duality}
With $\tau_0^*$ and $\nu_0$ defined on the residual rows,
\begin{equation}\label{eq:residual-duality}
\tau^*(N)=H(N)+\tau_0^*(N),\qquad
\nu(N)=H(N)+\nu_0(N).
\end{equation}
Moreover, for every $0<\theta<1/4$,
\begin{equation}\label{eq:fractional-lower}
\tau^*(N)\gg_\theta N^\theta
\end{equation}
for all sufficiently large even $N$.
\end{proposicion}
\begin{proof}
For a prime pair $\lambda<\varrho$, transfer the fractional weight at $\varrho$ to $\lambda$ and truncate the resulting weight at $1$. The row constraint ensures that this leaves weight $1$ at $\lambda$. No other row uses $\varrho$, so feasibility is preserved and cost cannot increase. The diagonal and endpoint force unit weights as before. This proves the first identity in \eqref{eq:residual-duality}.

A residual family of disjoint products can be joined to all the prime-pair rows and the endpoint. The primes in $\mathcal F_N$ do not occur in residual rows, and their larger prime partners occur only in their own rows. Conversely, any family of disjoint products has at most $H(N)$ rows using a prime from $\mathcal F_N$, since different rows cannot use the same such prime. Its other rows form a residual disjoint family. This proves the second identity.

For \eqref{eq:fractional-lower}, give each row of $\mathcal T_\theta(N)$ the weight
$\log^2N/(C_\theta\mathfrak S(N)N^{1-\theta})$ and every other row weight zero. Lemma~\ref{lem:prime-loss} makes these weights feasible for \eqref{eq:fractional-dual}. Summing them and using Lemma~\ref{lem:shifted-primes} proves the claim.
\end{proof}

\begin{ejemplo}[A fractional gap at $N=352$]\label{ex:fractional-gap}
Here $E_N=\{5+6j:0\le j\le28\}$, $\mathcal F_N=\{5,41,59,71,83,89,101,113,173\}$, and $H(N)=9$. Removing the rows in $\mathcal F_N$ and those with $5\mid h(N-h)$ leaves the following nine rows. Their displayed prime factors avoid $\mathcal F_N$, so the list is exactly $E_N^0$:
\begin{center}\small
\begin{tabular}{r l@{\qquad}r l}
\toprule
$h$&$\mathcal P_h$&$h$&$\mathcal P_h$\\\midrule
11&$\{11,31\}$&119&$\{7,17,233\}$\\
23&$\{7,23,47\}$&131&$\{13,17,131\}$\\
29&$\{17,19,29\}$&143&$\{11,13,19\}$\\
53&$\{13,23,53\}$&149&$\{7,29,149\}$\\
161&$\{7,23,191\}$&&\\\bottomrule
\end{tabular}
\end{center}
The primes $\{7,11,17,23\}$ cover these rows. A fractional cover gives weight $1$ to $11$ and weight $1/2$ to each of $7,13,17,23,29$, with total cost $7/2$. A dual solution of the same value gives weight $1$ to row $11$ and weight $1/2$ to each of rows $23,29,53,131,149$. The table verifies every constraint exactly. Thus $\tau_0^*=7/2$ and $\tau_0=4$. Rows $11,131,161$ have pairwise disjoint prime factors, so $\nu_0=3$, using $\nu_0\le\lfloor\tau_0^*\rfloor$. Therefore
\[
\nu(352)=12,\qquad\tau^*(352)=\frac{25}{2},\qquad\tau(352)=13.
\]
The displayed factorizations and rational weights prove a gap between integer covering and packing. No external data or numerical optimization is needed.
\end{ejemplo}

\subsection{A lower bound for the extremal size}
The cover reduction also constructs vanishing supports and therefore complements the upper bound for $J(k)$.

\begin{proposicion}\label{prop:J-lower}
As $k\to\infty$, one has $J(k)\gg k(\log k)^2$.
\end{proposicion}
\begin{proof}
For $N=2^j$, the local factor is $\mathfrak S(N)=1$. Apply Lemma~\ref{lem:upper-forms} with $q=1$ and a fixed $\theta$, for example $\theta=1/8$. Every prime pair whose smaller summand exceeds $w=N^{1/100}$ is counted by that upper sieve, apart from the possible diagonal. The pairs with smaller summand at most $w$, the diagonal, and the endpoint contribute $O(w+1)$. Hence
\[
H(N)\ll\frac{N}{\log^2N}.
\]
By Theorem~\ref{thm:cover-reduction},
\[
\tau(N)\le H(N)+M_0(N)\ll\frac{N}{\log^2N},
\]
since $M_0(N)\le\sqrt N=o(N/\log^2N)$.
Choose a sufficiently small fixed $\gamma>0$ and, for each large $k$, take the largest power of two $N\le\gamma k\log^2k$. Then $N>\gamma k\log^2k/2$, $\log N\sim\log k$, and $\tau(N)\le k$. A minimum cover yields $g_P(N)=0$. Every prime factor of $P$ is noncolliding because $N$ is a power of two. Thus $N\le J(k)$, proving the claim.
\end{proof}

Together with Corollary~\ref{cor:large-factors}, this proves Theorem~\ref{thm:main-J}. The lower bound constructs zeros for unrestricted supports. These covers need not satisfy $p_{\max}^2<N$.

\begin{observacion}[Forced noncolliding cost]\label{rem:forced-cost}
Every vanishing support satisfies the sharper inequality
\begin{equation}\label{eq:forced-cost}
k(N,P)\ge H(N)-\ind_{N/2\text{ prime}}\ge H(N)-1.
\end{equation}
Indeed, a prime pair $\lambda<\varrho$ with $\lambda+\varrho=N$ and $\lambda\ge5$ forces $P$ to contain $\lambda$ or $\varrho$. Neither divides $n$: one has $\varrho>n$, and $\lambda\mid n$ would imply $\lambda\mid\varrho$. Distinct pairs use distinct primes. The endpoint, when present, forces the noncolliding prime $N-1$. Only the diagonal prime $n$, when present, collides. Consequently,
\begin{equation}\label{eq:J-forced}
J(k)\le\sup\bigl(\{0\}\cup
\{N\ge26:N\text{ even},\ H(N)\le k+1\}\bigr).
\end{equation}
If $H(N)\gg N/\log^2N$ for all sufficiently large even $N$, as predicted by the Hardy--Littlewood conjecture \cite{KumchevTolev}, then \eqref{eq:forced-cost} gives $N\ll k(\log k)^2$ for every zero with at most $k$ noncolliding factors, for all sufficiently large $k$. Proposition~\ref{prop:J-lower} would then yield
\[
J(k)\asymp k(\log k)^2.
\]
Thus the conjectural exponent of $k$ is $1$. Proposition~\ref{prop:J-cutoff} uses exponent $2$ with a leading constant below $1$, while Theorem~\ref{thm:main-J} gives $4+\varepsilon$.
\end{observacion}
\section{The natural cutoff and the remaining gap}\label{sec:alcance}
The estimates for $J(k)$ quantify the size of a vanishing instance, while the cover reduction identifies its forced prime-pair component. Neither bound reaches the scale needed to settle positivity at the complete natural cutoff.

\begin{proposicion}[A sufficient quadratic bound]\label{prop:J-cutoff}
If there is a fixed $0<\eta<1$ such that
\begin{equation}\label{eq:J-quadratic}
J(k)\le(1-\eta)(k\log k)^2
\end{equation}
for all sufficiently large $k$, then $g_{\Pi(N)}(N)>0$ for every sufficiently large even $N$.
\end{proposicion}
\begin{proof}
Put $k_N=\omega(\Pi(N)/\gcd(\Pi(N),N/2))$. The prime number theorem and
$\omega(N)\ll\log N$ give
\[
k_N\sim\frac{2\sqrt N}{\log N},
\qquad (k_N\log k_N)^2\sim N.
\]
If $g_{\Pi(N)}(N)=0$, the definition of $J$ gives $N\le J(k_N)$, contradicting \eqref{eq:J-quadratic} for large $N$.
\end{proof}

The proposition is a sufficient criterion, not a converse or an established estimate for $J$. The exponent $4+\varepsilon$ in Theorem~\ref{thm:main-J} is separated from this quadratic target. The analogy with the ordinary Jacobsthal bound does not transfer its proof to the two classes tied to $N$.

The weighted criterion locates another precise gap. For the complete natural support, Mertens' theorem gives
\begin{equation}\label{eq:weight-natural}
W_\theta(N,\Pi(N))
=\sum_{\substack{N^\theta<q\le\sqrt N\\q\nmid N}}\frac1q
=\log\frac1{2\theta}+o_\theta(1).
\end{equation}
The factors dividing $N$ contribute only $O_\theta(N^{-\theta})$ in this range. Thus the upper bound for the total loss now has the same order as the supplied family, $\mathfrak S(N)N/\log^2N$. The constants obtained here do not prove
$c_\theta>C_\theta\log(1/(2\theta))$. Finite measurements of the actual losses cannot replace a uniform strict inequality. The gap is consistent with the parity obstruction in classical sieve theory. Our estimates give no lower bounds for intersections of losses. Type II, or bilinear, information is one possible route beyond these bounds, as illustrated by the prime-producing sieve framework in \cite{FordMaynard}. No such estimates are proved here.

For completeness, the local question remains:
\begin{problema}\label{prob:uniforme}
Does every even $N\ge26$, every prime $5\le p<\sqrt N$, and every product $Q$ of primes in $[5,p)$ satisfy $s_Q(N)>0$ and $\U(N,Q,p)$?
\end{problema}
An affirmative answer would exclude every first zero at the natural cutoff. It would also imply, by \eqref{eq:goldbach-endpoint}, that every even $N\ge26$ with $N-1$ composite is a sum of two primes exceeding $\sqrt N$. Because a large central occupancy can preserve a representation without satisfying $\U$, this local question contains more information than positivity alone.

Order the factors of a hypothetical vanishing support increasingly and let $P=Qp$ be its first vanishing prefix. The preceding set is nonempty. At the natural cutoff and for $N>2780460$, Theorems~\ref{thm:once} and~\ref{thm:siete} and Proposition~\ref{prop:bloques} require
\begin{equation}\label{eq:residual}
\omega(P)\ge12,\qquad k(N,P)\ge8,\qquad p^2<N<6Q,\qquad \ell<Q.
\end{equation}
The explicit results also require the profile inequalities \eqref{eq:factorial}--\eqref{eq:cardinal-n}.
For each fixed $\theta<1/4$, sufficiently large instances must have
$\kappa_\theta(N,P)\gg_\theta N^\theta$.
Corollary~\ref{cor:prefix-uniformity} excludes $N\ge p^{4+\varepsilon}$ for sufficiently large $p$, with any fixed $\varepsilon>0$.
The first-zero configuration is therefore constrained both by its residue concentration and by the cost of the support that produces it. Excluding every remaining instance is the open step.
\begingroup
\small
\setlength{\bibsep}{4pt}
\bibliographystyle{plainnat}
\bibliography{referencias}
\endgroup
\end{document}